\documentclass[a4paper,12pt]{article}

\usepackage{amsfonts,amssymb,amsmath,amsthm}

\newtheorem{theorem}{Theorem}

  \newtheorem{corollary}{Corollary}

 \def\dj {\hbox{\rlap{\kern 0.2em\raise 1.30ex\hbox
  {\vrule height 0.09ex width 0.33em}}d}}

\newcommand{\Z}{{\mathbb Z}}

\newcommand{\C}{{\mathbb C}}

\newcommand{\ep}{\varepsilon}

\newcommand{\er}{e_{r}^{*}}
\DeclareMathOperator{\ho}{H} 
\DeclareMathOperator{\sep}{sep}

\begin{document}

\title{{\Large{\bf Root separation for reducible monic polynomials of odd degree}}}

\author{Andrej Dujella and Tomislav Pejkovi\'{c}}

\date{{\it\small Dedicated to the memory of Professor Sibe Marde\v{s}i\'c}}
\maketitle

\begin{abstract}
We study root separation of reducible monic integer polynomials
of odd degree. Let $\ho(P)$ be the na\"\i ve height, $\sep(P)$ the minimal distance
between two distinct roots of an integer polynomial $P(x)$  and
$\sep(P)=\ho(P)^{-e(P)}$. Let $ \er(d)=\limsup_{\deg(P)=d,\, \ho(P)\to +\infty}e(P)$,
where the limsup is taken over the reducible monic integer polynomials $P(x)$ of degree $d$.
We prove that $\er(d) \leq d-2$. We also obtain a lower bound for $\er(d)$ for $d$ odd,
which improves previously known lower bounds for $\er(d)$ when $d\in\{5,7,9\}$.
\end{abstract}

\footnotetext{


The authors were supported by the Croatian Science Foundation under the project no.~6422. A.~D. acknowledges support from the QuantiXLie Center of Excellence.}


\section{Introduction}
All the polynomials that we deal with in this paper have integer coefficients. For any such polynomial,
we can look at how close two of its real or complex roots can be. Since we can always find polynomials
with distinct roots as close as desired, we need to introduce some measure of size for polynomials with
which we can compare this minimal separation of roots. This is done by bounding the degree and most
usually using the na\"\i ve height, that is, the maximum of the absolute values of the coefficients of a polynomial.

The problem of minimal root separation for polynomials with fixed degree has been completely solved
only in the trivial case of quadratic polynomials. Best possible separation exponent is also known
for nonmonic cubic polynomials (see \cite{Ev,Sch}). For monic cubic polynomials,
complete resolution would be equivalent to proving or disproving
the well known Hall's conjecture \cite{BuMi}. Therefore, resolving the problem completely for
polynomials of larger degree seems entirely out of reach.

However, if we restrict ourselves to reducible monic polynomials, then the cubic case becomes easy
and the quartic case has been solved by the authors \cite{DP11}. Thus, we are interested in
the separation properties of the reducible monic polynomials of degree at least $5$.

Let $P(x)$ be a polynomial of degree $d\geq 2$, na\"\i ve height $\ho(P)$ and with at least two distinct roots.
The \emph{polynomial root separation} of $P(x)$  is
\[ \sep(P):=\min_{\substack{P(\alpha)=P(\beta)=0,\\ \alpha\neq\beta}}|\alpha-\beta|. \]
The quantity $e(P)$ is defined by
\[ \sep(P)=\ho(P)^{-e(P)}. \]
Following the notation introduced in \cite{BuMi}, for $d\geq 2$, we set
\[ e(d):=\limsup_{\deg(P)=d,\, \ho(P)\to +\infty}e(P) \]
and
\[ \er(d):=\limsup_{\deg(P)=d,\, \ho(P)\to +\infty}e(P),  \]
where the latter limsup is taken over the reducible monic integer polynomials $P(x)$ of degree $d$.

There are some other variants of polynomial root separation problem like $p$-adic root separation,
where $|\alpha - \beta|$ is replaced by $|\alpha - \beta|_p$ (see \cite{Phazu}),
and absolute root separation, where $|\alpha - \beta|$ is replaced by $||\alpha| - |\beta||$ (see \cite{BDPS}),
but they will not be treated in this paper.

Obviously, we have $e(d)\geq \er(d)$. A classical result of Mahler \cite{Ma} says that $e(d)\leq d-1$ for every $d\geq 2$.
It is easy to see that $\er(2)=0 $ and $\er(3)=1$, while the main result of \cite{DP11} shows that
$\er(4)=2$. The best current lower bounds for the values we are interested in are
$\er(5)\geq 2$ from \cite{BuMi} and the following general result by Bugeaud and Dujella \cite{BD14}
\[
\er(d)\geq \frac{2d}{3}-1 \textrm{ for even } d\geq 6, \quad
\er(d)\geq \frac{2d}{3}-\frac{5}{3} \textrm{ for odd } d\geq 7.
\]
In particular, this implies that $\er(7)\geq 3$ and $\er(9)\geq \frac{13}{3}$.
The mentioned result from \cite{BD14} is obtained by constructing the parametric family
of polynomials $Q_{d,n}(x)=(x^2-(n^2+3n+1)x+(n+2))q_{d-2,n}(x)$, where $q_{d-2,n}(x)$ is a recursive sequence
of polynomials of even degree. For $d$ odd, the polynomials
$Q_{d,n}(x)=x(x^2-(n^2+3n+1)x+(n+2))q_{d-3,n}(x)$ were used, so it is not surprising that results
for odd degrees are sightly weaker and, in particular, for $d=5$ the lower bound $\er(5)\geq \frac{5}{3}$
from \cite{BD14} is weaker than the bound $\er(5)\geq 2$ from \cite{BuMi}.

In this paper, we mainly consider odd degree polynomials.
In Section \ref{sec:2} we improve the upper bound $\er(d)\leq d-1$ which follows from \cite{Ma}
by proving that $\er(d)\leq d-2$ for $d\geq 2$.
In Section \ref{sec:3}, we will construct a parametric family of reducible monic polynomials
of odd degree with good root separation properties. Although the obtained lower bound will be
asymptotically weaker compared with the bound from \cite{BD14}, it will improve all previously
known lower bounds on $\er(d)$ for $d\in\{5,7,9\}$. We show that
\[
  \er(5) \geq \frac{7}{3}, \quad
  \er(7) \geq \frac{17}{5}, \quad
  \er(9) \geq \frac{31}{7}.
\]

\section{Upper bounds}  \label{sec:2}

Our first goal is to we improve the upper bound $\er(d)\leq d-1$ which follows from \cite{Ma}
by proving that $\er(d)\leq d-2$ for $d\geq 2$. Note that we already know that $\er(d)= d-2$ for $d=2,3,4$,
so we may assume $d\geq 5$ here.

Let $P(x)$ be a monic polynomial of degree $d$, reducible over the rationals.
Gauss's Lemma shows that we can factor $P(x)$ into two nonconstant monic polynomials with integer coefficients.
If $P(x)=Q(x)R(x)$, where $Q(x)$ and $R(x)$ are integer polynomials of positive degrees $n$ and $m$, respectively, and $\alpha,\beta\in\C$ are such that $Q(\alpha)=R(\beta)=0$, then a version of Liouville's inequality (see \cite[Theorem A.1]{Bu}) ensures that
\begin{equation} \label{eq:liouville}
  |\alpha - \beta| \gg \ho(Q)^{-m} \ho(R)^{-n},
\end{equation}
where the constant implied in the Vinogradov symbol $\gg$ here and further on depends only on $m$ and $n$.

Gelfond's Lemma \cite[Lemma A.3]{Bu} says that
\begin{equation} \label{eq:gelfond}
  2^{-m-n}\ho(Q)\ho(R) \leq \ho(P) \leq 2^{m+n} \ho(Q)\ho(R).
\end{equation}
Now \eqref{eq:liouville} and \eqref{eq:gelfond} and the fact that $\ho(Q)\geq 1$, $\ho(R)\geq 1$ give
\[ |\alpha - \beta| \gg \ho(P)^{-\max\{m,n\}}. \]
Mahler's general result shows that
\begin{equation} \label{eq:Q}
  \sep(Q) \gg \ho(Q)^{-n+1} \gg \ho(P)^{-n+1}
\end{equation}
and
\begin{equation} \label{eq:R}
  \sep(R) \gg \ho(R)^{-m+1} \gg \ho(P)^{-m+1}.
\end{equation}
Combining the last three inequalities, we obtain
\begin{equation} \label{eq:P}
  \sep(P) \gg \ho(P)^{-\max\{m,n\}}.
\end{equation}

If $\min\{m,n\}\geq 2$, we have $\max\{m,n\}\leq d-2$ and \eqref{eq:P} shows that $\sep(P)\gg \ho(P)^{-d+2}$.

We see that the only case that needs to be considered more thoroughly in order to prove the upper bound
$\er(d)\leq d-2$ is when one of the polynomials $Q(x)$, $R(x)$ is linear and we have a root of $Q(x)$
and a root of $R(x)$ which are very close. Without loss of generality, assume that $Q(x)=x-c$,
$R(c)\neq 0$ and $R(c+\ep)=0$, where $c$ is an integer and $\ep$ is small, but nonzero, say $0<|\ep|<1$.

An easy bound $|c+\ep|\leq m\ho(R)$ is obtained by substituting $x=c+\ep$ into $R(x)$ and
comparing the leading term with the rest.

If $|c+\ep|\geq 3$, then $|c+\ep|\geq |c|-|\ep|\geq |c|-1\geq |c|/2$ since
$|c|\geq |c+\ep|-|\ep|\geq |c+\ep|-1\geq 2$. The first inequality in \eqref{eq:gelfond} together with $\ho(Q)=\max\{1,|c|\}$ and $|c|/2\leq |c+\ep| \leq m\ho(R)$ gives
\begin{equation} \label{eq:sqrtHP}
  |c|\ll \ho(P)^{1/2}.
\end{equation}
In case $|c+\ep|<3$, we have $|c|\leq |c+\ep|+|\ep|<4$ and since $\ho(P)\geq 1$, inequality \eqref{eq:sqrtHP} also holds.

Rolle's mean value theorem gives
\begin{equation} \label{eq:eps}
  1\leq |R(c)| = |R(c+\ep)-R(c)| = |R'(t)|\cdot |\ep|,
\end{equation}
where $t$ is between $c$ and $c+\ep$ and thus in the interval $(c-1,c+1)$. Since
\[ \ho(R')\leq  m\ho(R) \ll \ho(P) \quad\textrm{and}\quad
  |t|<|c|+1\ll \ho(P)^{1/2}, \]
we have
\[ |R'(t)|\ll \ho(P)^{\frac{m-1}{2}} + \ho(P)\cdot\ho(P)^{\frac{m-2}{2}} \ll \ho(P)^{\frac{m}{2}}, \]
where we used the fact that $R(x)$ is monic and $R'(x)$ is of degree $m-1$. Comparing with \eqref{eq:eps}, we obtain
\begin{equation} \label{eq:linear}
  |\ep|\gg \ho(P)^{-m/2} \gg \ho(P)^{-(d-1)/2} \gg \ho(P)^{-d+3}
\end{equation}
since $d\geq 5$ (we get $|\ep|\gg \ho(P)^{-(d-1)/2} \gg \ho(P)^{-d+2}$ for $d\geq 3$).

Hence, we proved the following theorem.

\begin{theorem} \label{th:2}
For $d\geq 2$, it holds that $\er(d)\leq d-2$.
\end{theorem}

\section{Lower bounds} \label{sec:3}
We would like to improve the lower bound $\er(5)\geq 2$ from \cite{BuMi}.
By \eqref{eq:linear}, we see that this bound cannot be improved by
considering degree $5$ polynomials which have a linear factor.
Thus, we have to consider products of two monic polynomials of degrees $2$ and $3$.

We made some experiments in order to find suitable degree $5$ polynomials with $e(P)>2$.
We consider monic cubic polynomials with coefficients of moderate size. For such a polynomial, we
choose one of its real roots $\gamma$
and then apply the LLL-algorithm to a matrix of the form
$$ \left(
  \begin{array}{ccc}
    1 & 0 & 0 \\
    0 & 1 & 0 \\
    N & \lfloor N \gamma \rceil & \lfloor N \gamma^2 \rceil \\
  \end{array}
\right),
$$
with suitably chosen large positive integer $N$,
to construct quadratic polynomials with a root close to $\gamma$
(the method is explained e.g. in \cite[Chapter 6]{Smart}).
Among the polynomials obtained in the experiments,
we have noted two collections of polynomials with $e(P)$ approaching $\frac{7}{3}$:
$$ s_n(x)=(x^3-2nx^2+(2-2n)x+2)(x^2+(-2n^2-2n)x+2n+2) $$
and
$$ p_n(x)=(x^3+nx-1)(x^2+n^2x-n). $$
Indeed, $s_n(x)$ has two close roots with asymptotic expansions
$$ \frac{1}{n}+\frac{1}{2n^4} -\frac{1}{2n^5}+ \frac{1}{2n^6} + O(\frac{1}{n^8}) \quad
\mbox{and} \quad  \frac{1}{n}+\frac{1}{2n^4} -\frac{1}{2n^5}+ \frac{1}{2n^6} + \frac{1}{4n^7} + O(\frac{1}{n^8}), $$  while $p_n(x)$ has two close roots with asymptotic expansions
$$ \frac{1}{n}-\frac{1}{n^4}+\frac{3}{n^7} + O(\frac{1}{n^{10}}) \quad
\mbox{and} \quad \frac{1}{n}-\frac{1}{n^4}+\frac{2}{n^7} + O(\frac{1}{n^{10}}). $$
Hence, we proved that $\er(5) \geq \frac{7}{3}$.
Moreover, it is easy to see that the polynomials $p_n(x)$
can be generalized to arbitrary odd degree.

Let $d\geq 5$ be an odd number. To obtain a lower bound on $\er(d)$,
we construct a family $(P_{d,n})_{n\geq 1}$
of reducible monic polynomials of degree $d$ (such that $P_{5,n}=p_n$),
depending on the parameter $n$,
with root separation asymptotically
\begin{equation} \label{eq:sepP}
  \sep(P_{d,n})\ll \ho(P_{d,n})^{-\frac{d^2-2d-1}{2d-4}},\quad n\to +\infty.
\end{equation}
This will give
\begin{equation} \label{eq:sepPb}
\er(d)\geq\frac{d^2-2d-1}{2d-4} = \frac{d}{2} - \frac{1}{2d-4}.
\end{equation}
The bound \eqref{eq:sepPb} is comparable with best known lower bounds
for separation of irreducible monic
and nonmonic polynomials ($e_{irr}^{*}(d) \geq \frac{d}{2} - \frac{1}{4}$,
$e_{irr}(d) \geq \frac{d}{2} + \frac{d-2}{4(d-1)}$ see \cite{BD14,BD11}).
Although it is asymptotically weaker than the bound $\er(d)\geq (2d-5)/3$ from \cite{BD14},
it is better than $\er(d)\geq (2d-5)/3$ for $d=5$, $d=7$ and $d=9$,
and it is also better than $\er(5)\geq 2$ from \cite{BuMi}.

For $k,n\in\Z_{\geq 2}$, let
\[ \begin{aligned}
   Q(x) &= Q_{k,n}(x) = x^2 + n^k x - n, \\
   R(x) &= R_{k,n}(x) = \frac{x^{2k+1}-Q(x)}{x^2-n} = x\cdot \frac{x^{2k}-n^k}{x^2-n} - 1 \\
     &= x(x^{2(k-1)} + nx^{2(k-2)} + \cdots + n^{k-1}) - 1.
\end{aligned} \]
We omit $k$ and $n$ from the index of polynomials $Q_{k,n}(x)$ and $R_{k,n}(x)$ for easier writing. 
It is clear that $Q(x)$ and $R(x)$ are monic integer polynomials of degrees $2$ and $2k-1$ and heights $n^k$ and $n^{k-1}$, respectively. Thus $Q(x)R(x)$ is a reducible monic polynomial of degree $2k+1$ and height $n^{2k-1}-1$ which is attained only in the quadratic term.

Quadratic polynomial $Q(x)$ has roots $(-n^k\pm \sqrt{n^{2k} + 4n})/2$, so that one root is close to $-n^k$ and the other root, which we denote by
\[ \alpha = \frac{2n}{n^k+\sqrt{n^{2k}+4n}}, \]
is close to $n^{1-k}$, more precisely, $0 < \frac{1}{2}n^{1-k} < \alpha < n^{1-k}.$

We also have
\[ R(\alpha) = \frac{\alpha^{2k+1}}{\alpha^2 - n} <0, \quad
  R(2\alpha) = \frac{(2\alpha)^{2k+1} - 2\alpha^2 - n}{(2\alpha)^2 - n} >0 \]
since the terms other than $-n$ in the numerator and the denominator are very small.
Therefore, the polynomial $R(x)$ has a root $\beta$ in the interval $(\alpha,2\alpha)$.

Using the mean value theorem, we conclude that there is $t\in(\alpha,\beta)$ such that
\begin{equation} \label{eq:rolle}
   0 < R'(t) (\beta-\alpha) = R(\beta) - R(\alpha) = \frac{\alpha^{2k+1}}{n - \alpha^2}
     < 2n^{(2k+1)(1-k)-1}.
\end{equation}
It is easily seen that
$R'(t)>n^{k-1}$ and employing this inequality in \eqref{eq:rolle} gives
\[ 0 < \beta - \alpha < \frac{2n^{(2k+1)(1-k)-1}}{n^{k-1}} = 2n^{1-2k^2}. \]
This implies
\[ \sep(Q\cdot R) < 2\ho(Q\cdot R)^{\frac{1-2k^2}{2k-1}}.  \]

For an odd integer $d\geq 5$, we take $k=(d-1)/2$ and define $P_{d,n}(x):=Q_{k,n}(x)R_{k,n}(x)$ obtaining a family
that satisfies \eqref{eq:sepP}.

\begin{corollary} \label{kor:3}
It holds that
$$  \frac{7}{3} \leq \er(5) \leq 3, \quad
  \frac{17}{5} \leq \er(7) \leq 5, \quad
  \frac{31}{7} \leq \er(9) \leq 7. $$
\end{corollary}

\bigskip


\bigskip

\bigskip

\begin{center}
{\bf Separacija korijena za reducibilne normirane polinome neparnog stupnja}
\end{center}

\bigskip

\begin{center}
{\it Andrej Dujella i Tomislav Pejkovi\'c}
\end{center}

\bigskip

\begin{center}
\begin{minipage}[c]{9.2cm}
{\small \hspace*{0.5cm} {\sc Sa\v{z}etak.}
U ovom \v{c}lanku prou\v{c}avamo separaciju korijena reducibilnih normiranih
polinoma neparnog stupnja. Neka je $\ho(P)$ visina, $\sep(P)$ minimalna uda\-lje\-nost
razli\v{c}itih korijena polinoma $P(x)$ s cjelo\-broj\-nim koeficijentima, te
$\sep(P)=\ho(P)^{-e(P)}$. Neka je $ \er(d)=\limsup_{\deg(P)=d,\, \ho(P)\to +\infty}e(P)$,
gdje se limsup uzima po svim reducibilnim normiranim cjelobrojnim polinomima $P(x)$ stupnja $d$.
Dokazujemo da vrijedi $\er(d) \leq d-2$. Tako\dj{}er, dobivamo donju ogradu za $\er(d)$ za neparan $d$,
koja poboljs\v{s}ava prethodno poznate donje ograde za $\er(d)$ za $d\in\{5,7,9\}$.}%
\end{minipage}
\end{center}

\bigskip

{\small \noindent
Andrej Dujella \\
Department of Mathematics \\
Faculty of Science\\
University of
Zagreb
\\ Bijeni\v{c}ka cesta 30 \\
10000 Zagreb, Croatia \\
{\em E-mail address}: {\tt duje@math.hr}}

\bigskip

{\small \noindent
Tomislav Pejkovi\'c \\
Department of Mathematics \\
Faculty of Science \\
University of
Zagreb
\\ Bijeni\v{c}ka cesta 30 \\
10000 Zagreb, Croatia \\
{\em E-mail address}: {\tt pejkovic@math.hr}}

\end{document}